\documentclass{article}

\usepackage{amsfonts,xr,graphicx,amsmath}

\begin{document}
\title{{Homogeneous  Solutions of Fully Nonlinear Elliptic
 Equations in Four Dimensions} }
 \author{{Nikolai Nadirashvili\thanks{LATP, CMI, 39, rue F. Joliot-Curie, 13453
Marseille  FRANCE, nicolas@cmi.univ-mrs.fr},\hskip .4 cm Serge
Vl\u adu\c t\thanks{IML, Luminy, case 907, 13288 Marseille Cedex
FRANCE, vladut@iml.univ-mrs.fr} }}

\date{}
\maketitle

\def\S{\mathbb{S}}
\def\Z{\mathbb{Z}}
\def\R{\mathbb{R}}
\def\N{\mathbb{N}}
\def\H{\mathbb{H}}
\def\tilde{\widetilde}
\def\epsilon{\varepsilon}

\def\n{\hfill\break} \def\al{\alpha} \def\be{\beta} \def\ga{\gamma} \def\Ga{\Gamma}
\def\om{\omega} \def\Om{\Omega} \def\ka{\kappa} \def\lm{\lambda} \def\Lm{\Lambda}
\def\dl{\delta} \def\Dl{\Delta} \def\vph{\varphi} \def\vep{\varepsilon} \def\th{\theta}
\def\Th{\Theta} \def\vth{\vartheta} \def\sg{\sigma} \def\Sg{\Sigma}
\def\bendproof{$\hfill \blacksquare$} \def\wendproof{$\hfill \square$}
\def\holim{\mathop{\rm holim}} \def\span{{\rm span}} \def\mod{{\rm mod}}
\def\rank{{\rm rank}} \def\bsl{{\backslash}}
\def\il{\int\limits} \def\pt{{\partial}} \def\lra{{\longrightarrow}}

 {\em Abstract.} We prove that there is no nontrivial homogeneous order 2 solutions
 of fully nonlinear uniformly  elliptic equations in dimension 4.

\section{Introduction}
\bigskip

We study  a class of solutions to fully
nonlinear second-order elliptic equations of the form
\begin{equation}F(D^2u)=0\end{equation}
$D^2u$ being the Hessian of the function $u$ defined in  $\R^n$. We assume that
$F$ is a smooth function defined on the space $ S^2(\R^n)$
of ${n\times n}$ symmetric matrices  satisfying the uniform ellipticity condition:
$${1\over C'}|\xi|^2\le F_{u_{ij}}\xi_i\xi_j\le C' |\xi |^2\;,
\forall\xi\in {\bf R}^n\;.$$
 Here, $u_{ij}$ denotes the partial derivative
$\pt^2 u/\pt x_i\pt x_j$. A function $u$ is called a {\it
classical\/} solution of (1) if $u\in C^2(\Om)$ and $u$ satisfies
(1).  Actually, any classical solution of (1) is a smooth
($C^{\alpha +3}$) solution, provided that $F$ is a smooth
$(C^\alpha )$ function of its arguments.

Let $B=\{ x\in \R^n: |x|<1 \}$ be a ball, $g$ be a continuous function on $\partial B$.
Consider a Dirichlet problem

\begin{equation}\left\{
\begin{array}{l l}
 F(D^2u)=0 &\mbox{in } B \\
u=g &\mbox{on $\partial B\;$} \\ \end{array} \right. \end{equation}

We are interested in the problem of existence and regularity of
solutions to the Dirichlet problem (2).
The problem (2) has always a unique viscosity (weak)
solution for  fully nonlinear elliptic equations. The viscosity solutions  satisfy the equation
(1) in a weak sense, and the best known interior regularity
([C],[CC],[T3]) for them is $C^{1,\epsilon }$ for some $\epsilon>
 0$. For more details see [CC], [CIL]. Note, however, that viscosity solutions
are $C^{2,\epsilon }$-regular  almost everywhere; in fact, it is true on the complement
of  a closed set of Hausdorff dimension strictly less then $n$ [ASS].
Until recently it remained unclear whether non-smooth viscosity solutions exist. In the recent papers
 [NV1], [NV2],  [NV3], [NV4]  the authors first proved   the existence   of non-classical viscosity solutions
 to a fully nonlinear elliptic equation, and  of singular solutions to Hessian (i.e. dependinding only on the 
eigenvalues of $D^2u$) uniformly elliptic equation in all dimensions beginning from 12, and, finally, 
the paper   [NTV] gives a construction  of  non-smooth viscosity solution in 5 dimensions which is order 2
 homogeneous, also for Hessian equations. These papers   use the functions
$$w_5(x)={P_{5} (x)\over |x|},\; w_{12,\delta}(x)={P_{12}(x)\over |x|^{\delta }},\;w_{24,\delta}(x)= {P_{24} (x)\over  |x|^{\delta }},\:
\delta\in [1,2[,
$$
for certain (minimal) cubic forms  $P_{5} (x), P_{12}(x),P_{24}(x)$  in the dimensions 5,12 and 24, respectively.

\smallskip

On the other hand the classical Alexandrov's theorem [A] says  that an analytic in $\R^3\setminus \{ 0 \}$ homogeneous order 1 function $u$ such that the Hessian $D^2u$ is either non-definite or
0 at any point is linear. This immediately implies the absense of  homogeneous order 2 real analytic 
 in $\R^3\setminus \{ 0 \}$ solutions to fully nonlinear equations  different from quadratic forms
 (in   $C^{2, \alpha}$ setting it is proved in [HNY]).  Thus  the existence of  homogeneous order 2 real
 analytic  outside zero solutions to fully nonlinear equations is not known exactly in 4 dimensions,
the analogue of Alexandrov's theorem in 4 dimensions being false (indeed 
$u = (x_1^ 2+x_2^ 2-x_3^ 2-x_4^ 2)/|x|$ gives a counter-example, cf. [LO]).\smallskip\smallskip
  
This note   fills this gap showing that    5  is the minimal dimension where there exist homogeneous
 order 2 non-smooth  solutions to uniformly elliptic fully nonlinear equations.   

\medskip
{\bf Theorem 1.}  {\it Let $u$ be a homogeneous order $2$ real analytic function in $\R^4 \setminus \{0\} $.
 If $u$ is a solution of the uniformly  elliptic equation $F(D^2u)=0 $ in $\R^4 \setminus \{0\} $,
then $u$ is a quadratic polynomial. }

\medskip We collect some preliminary lemmas in Section 2 below and give the proof in Section 3.

\section{Preliminary results}

Here we prove some general results we need to prove the theorem. 

\medskip
{\bf Lemma 0. } {\it Let $v$ be a smooth homogeneous order $1$ function in $\R^3\setminus \{ 0 \}$.
Assume that $y\in  \S^2 $ and the quadratic form $D^2v(y)$ changes sign. Let $a\in \S^3 , a\neq y$, 
and let $G$ be an open domain in $\R^3,y\in G$. Then
$$\sup_G v_a(x) >v_a(y).$$ }

\medskip
{\em Proof.} Let $L\subset \R^3 $ be an affine 2-dimensional plane transversal to the vector $y$ such
that $y\in L$ and $a$ is parallel to $L$. Denote by $v'$ the restriction of the function $v$ on $L$. Since
$v$ is a homogeneous order 1 function the quadratic form $D^2v'(y)$ changes sign. Thus there
is a neighborhood  $D$ of the point $y$ where $v'$ satisfies a uniformly elliptic equation on $L$
of the form
 $$ \sum a_{ij} (x){\partial^2v' \over  \partial x_i \partial x_j } =0.$$
 Thus by the maximum principle for the gradient of a solution of elliptic equations in dimension 2, 
see [GT], $v'_a$ cannot attain the supremum at the point $y$. The lemma is proved.

\medskip{\bf Lemma 1. } {\it Let $v$ be a real analytic homogeneous order $1$ function in
 $\R^n\setminus \{ 0 \}$. Assume that $v$ is a solution of a linear uniformly elliptic equation
 $$Pv=  \sum a_{ij} (x/ |x|){\partial^2v \over  \partial x_i \partial x_j } =0,$$
where coefficients $a_{ij}$ are smooth functions on $\S^{n-1}$. Let $ e_1,...,e_n\in \S^{n-1}$  be
 linearly independent unit vectors. Assume that the functions $v_{e_i}$, $i=1,...,n$ attain local 
supremum at $a\in \S^{n-1},  a\neq e_i, i=1,..., n$. Then $v$ is a linear function. }

\medskip{\em Proof.} Denote by $L$ an affine hyperplane in $\R^n$ orthogonal to $a$,
$ a \in L$. Then the restriction $v'$ of the function $v$ on $L$ satisfies a linear uniformly elliptic
 equation of the type
 $$P(v')= \sum a'_{ij} (y){\partial^2v' \over  \partial x_i \partial x_j } =0,$$
 where $y\in L$ and $a'_{ij}$ are smooth functions on $L$. Indeed, $D^2v(a)=0$  since $v$ is order one
 homogeneous, thus the partial derivatives  of $v'$ coinside with ones of $v$ in an appropriate
coordinate system. We consider then a coordinate system  on   $L$ such that the point $ a$ becomes
 the origin, assuming without loss
 that $v'(0)=0, \nabla v'(0)=0$. After a linear transformation of $\R^n$ we can assume  
 that   $P(0)$  is the Laplacian, i.e., $a'_{ij}(0)=\delta_i^j$.
Let  $p, \:\deg p=k\ge 2$ be the first nonzero homogeneous polynomial of the Taylor expansion of 
$v'$ at 0; clearly $p$ is  harmonic. Let $B\subset L$ be a  small ball centered at $0,$ 
let $g$ be the gradient map 
$$g:L\rightarrow \R^{n-1},\;g:= \nabla v'$$
 and let $\Gamma = g(B)$. Then $\Gamma \subset K:= \bigcap_{i=1}^n\{ e_i\leq 0\}$, $K$ being a strictly
 convex cone in $ \R^n$ since $e_i$ are linearly independent. Denote $K_0 =\{ K+a\}\cap L $;  if $K_0$ is 
non-empty then $K_0$ is a strictly convex cone in $L$. Let $p'$ be a non-zero partial derivative 
 of $p$ of order $k-2$; the quadratic form $p'$
changes sign, hence $\nabla p' (L)$ intersects the complement of $K_0$ and thus
  $ l^+\bigcap K_0=\emptyset$ for a line $l\subset \nabla p' (L)$ and 
a ray $ l^+\subset l$. Let $\Lambda:= \nabla {p'}^{-1}(l^+ ) $, then the curve
 $g(\Lambda ) \subset \R^n$ is
 tangent to  $l^+$ at the point $\{ a\}$ since $v_a (x)=O(|a-x|^{k})$.   Therefore
$g( \Lambda \cap B)$ intersects the complement of $K$, and the lemma follows.

\medskip{\bf Lemma 2. } {\it Let $v$ be a real analytic homogeneous order $1$ function in $\R^4\setminus \{ 0 \}$.
Assume that $v$ is a solution of a linear uniformly elliptic equation
\begin{equation}Pv=  \sum a_{ij} (x/ |x|){\partial^2v \over  \partial x_i \partial x_j } =0, \end{equation}
and the rank of the gradient map $\nabla v: \S^3 \rightarrow \R^4$ is $\leq 2$. Then
$v$ is a linear function. } 

\medskip{\it Proof.} Let $y\in \S^3$, $m\subset \R^4$ be a subspace,   $m\perp y$. Let $M\subset R^4$
be an affine hyperplane parallel to $m$, $y\in M$, and let  $f$ be the restriction of $v$ on $M$. 
Then $f$ is a real analytic function on $M$ such that for any $x\in M$ the hessian
$D^2f(x)$ is degenerate and either the quadratic form $D^2f(x)$ changes sign or $D^2f(x)=0$.
Let $$H:=\{x\in \R^3:\rank (D^2f(x))=2\}.$$ 
 We assume without loss that  $codim(\R^3 \setminus H)\ge 1.$  
 For $x\in H$ let $z(x)$ be the zero eigenspace of $D^2f(x)$. By assumption of the lemma $z(x)$ is a 
line analytically  depending  on the point $x\in H$.
 By Chern-Lashof's lemma,
 [CL, Lemma 2], \hskip .5cm
[S, Lemma VI 5.1] in the neighborhood of any point $x\in M$ 
 the plane $M$ is foliated by a 2-dimensional family of straight lines $L$, such that for any  line $l\in L$
  the restriction of the function $f$ on $l$ is an affine function, moreover  $l$ is parallel to
 the line $z(x)$ at any point $x\in l$, see the proof of Lemma 2 in [CL]. By the analyticity of $f$ it
 follows that the family $L$ foliate the whole space $M$ without intersection. 
  Let $l\in L$ and $p\subset \R^4$ be a two-dimensional plane spanned by  $l$ in
 $ \R^4$. Since $v$ is a homogeneous order one function it follows that  $v$ is linear on a half-plane of $p$.
 By analyticity, $v$ is a linear function on the whole plane $p$. Denote 
the whole set of these planes $p$ by $P$. Then any two planes of $P$
 intersect only at $\{ 0\}$ and foliate $\R^4\setminus m$. 
 
 Let $y'\in \S^3,\; m'=(y')^{\perp}\subset \R^4$ and let $P'$ be the foliation of $\R^4\setminus m'$ by  two-dimensional planes corresponding to $y'$. We will prove that $P$ and $P'$ coincide  on 
 $\R^4\setminus (m\cup m')$. Assume not. Then there is a 4-dimensional subset $X\subset  \R^4$
 such that for any $x\in X$ one has $x\in p\cap p'$ for some $p\in P$, $p'\in P'$, $p\neq p'$. 
Since the  planes $p$ and $p'$ are zero eigenspaces  of $D^2v$ it follows that
 that the zero eigenvalue  has multiplicity at least 3 at $x$, and hence $D^2v(x)=0$. Thus $D^2v$
 vanishes on $X$ and hence by analyticity of $v$ it follows that $v$ is a linear function. 
  Thus choosing different $y\in \S^3$ we get a foliation $P$ of $\R^4\setminus \{ 0\}$ by two
  dimensional planes which are zero eigenspaces of $D^2v$. 
  
 Notice that any 3-dimensional subspace of $\R^4$ contains at most one plane of $P$, since
 any two different planes in 3-dimensional space have nontrivial intersection. 
 
 Let $m\in \R^4$ be a 3-dimensional subspace such that $m\supset p,\: p\in P$. Denote by $v'$
 the restriction of the function of $v$ to $m$; subtracting a linear function  we can assume
 that $v'=0 $ on $p$. Let $x\in m\setminus p, \;x\in p'$ for some $p'\in P$. Then $p'$ is trasversal to $m$.
 Since $p'$ is a zero eigenspace of $D^2v$ it follows that either $D^2v' (x)$ changes sign or $D^2v'(x)=0$ on
 $m$. Thus the function $v'$ is a solution of an elliptic equation (3) at $x$. Thus we proved 
 that $v'$ satisfies an elliptic equation (3) on $m\setminus p$. Let  $e\in m$ be a vector parallel 
 to $p$. Let $z\in \S^2\subset m$ be a point at which   $v'_e $ attains its maximum on $\S^2$. 
If $v'_e(z)>0$, then $z\in \S^2\setminus p$ since by our assumption $v'=0$ on $p$. 
Since in a neighborhood of $z$ the function $v'$ is a solution of  (3) this contradicts  Lemma 0. Thus
 $v'_e(z)\le 0$ and thus  $v'_e\le 0$ everywhere since $v'_e(z)$ is maximal. Applying the same 
argument to the function $-v'$ we get $v'_e\ge 0$ everywhere and thus $v'_e \equiv 0$ for
 any vector $e$ parallel to $p$. Hence $v'$ is a function which depends only on the
 coordinate orthogonal to $p$ and therefore $v'$ is a linear function. Thus we get
 that for any three dimensional subspace $m$ of $\R^4$ the restriction of $v$ on $m$
 is a linear function. Hence $v$ is a linear function on $\R^4$ and the lemma is proved.

\medskip
{\bf Lemma 3.} {\it Let $Q(x,y,z)\in  \R[x,y,z]$ be  a cubic form such that for any $e\in {\S}^2 $ the quadratic
form $Q_e$ is degenerate. Then $Q$ is a function of two variables in some coordinate system.}

\medskip{\em Proof.} First of all, the conditions as well as the conclusion of the lemma are invariant under
non-singular linear transformations.  Considering $Q(x,y,z)=0$ as an equation of a plane projective
cubic curve $E_Q$  and applying the usual argument giving its Weierstrrass form
 (see, e.g. pp. 45-46 in the proof of Proposition 1.2 of Ch. 2 in [M]) one gets one the following:

\smallskip 1.  $E_Q$ is elliptic or irreducible possesing a singular point with $y\neq 0$;
in this case $Q$ is equivalent under a linear transfomation to the Weierstrass form
$$Q_W=y^2z+x^3+px^2z+qz^3;$$

\smallskip 2. $E_Q$ is  irreducible possesing a singular point with $y= 0$; then
$$Q=Q_s=x^3+axyz+bxz^2+cyz^2+dz^3$$
after a suitable non-singular linear transformation;

\smallskip
3. $E_Q$  is reducible, then either
$$Q=Q_r= z(x^2+ay^2+bz^2+cxz+dyz)$$
modulo such a transformation or $Q$ verifies the conclusion.

\smallskip \smallskip If $Q=y^2z+x^3+px^2z+qz^3, e=(k,l,m)$ then
$$r=r(k,l,m):=\det(D^2(Q_e))=3k^2mp+9km^2q-3kl^2-m^3p^2$$
should be indentically zero; in particular, $-6=r_{kll}=\partial^3 r/\partial k\partial l^2=0 $
which is clearly not the case.
\smallskip

If $Q=x^3+axyz+bxz^2+cyz^2+dz^3$ then
$$r/2=a^3klm+a^2bkm^2+a^2clm^2-3a^2dm^3+4abcm^3-3a^2k^3-12ack^2m-12c^2km^2,$$
$$0=r_{klm}/2=a^3,\;0=r_{kmm}/4=a^2b-12c^2$$
implying $c=a=0$ and the conclusion.

\smallskip If $Q= z(x^2+ay^2+bz^2+cxz+dyz)$ then
$$r/8=3abm^3-ac^2m^3-a^2l^2m-ackm^2-adlm^2-d^2m^3-ak^2m,$$
$$ 0=r_{llm}/16= -a^2,\;0=r_{mmm}/48= 3ab-ac^2-d^2$$
thus $a=d=0$ as necessary and the proof is finished.

\medskip
{\bf Lemma 4.} {\it Let $Q(x,y,z)\in  \R[x,y,z]$ be  a cubic form such that for any $a\neq b \in C\subset {\S}^2$ 
  the partial derivative  $Q_{ab}$ vanishes as  a linear form, $C$ being a curve on $ {\S}^2$ .
 Then $Q$ is a function of two variables in some coordinate system.}

\medskip{\em Proof.} The proof is very similar to that of Lemma , but sightly more combersome. We consider the 
same three main cases, each of them being divided in subcases depending on the curve $ C\subset {\S}^2$.

\smallskip
1). Weierstrass case. There are two subcases:

\smallskip \hskip .5 cm 1a). The curve $C$  is not in ${\S}^2\bigcap (\{y=0\}\bigcup\{z=0\}).$

\smallskip
\hskip .5 cm 1b). The curve $C\subset {\S}^2\bigcap (\{y=0\}\bigcup\{z=0\}).$

\smallskip
 In the subcase 1a we can suppose without loss that $a=(a_1,b_1,c_1),$\linebreak $ b=(a_2,b_2,c_2)$ with 
$c_1b_2+c_2b_1\neq 0$. A brute force calculation gives $Q_{aby}/2=c_1b_2+c_2b_1\neq 0$ and thus we get a contradiction.

\smallskip
 In the subcase 1b we  suppose without loss that $a=(a_1,b_1,0),\: b=(a_2,b_2,0)$ with $a_1a_2\neq 0$
but then $Q_{abx}/6=a_1a_2 \neq 0.$

\smallskip\smallskip
2). Singular case (singularity at $y=0$), $Q=x^3+pxyz+qxz^2+ryz^2+sz^3$. Subcases:

\smallskip \hskip .5 cm 2a). The curve $C$  is not in ${\S}^2\bigcap\{z=0\}.$

\smallskip
\hskip .5 cm 2b). The curve $C\subset {\S}^2\bigcap\{z=0\}.$

\smallskip Suppose 2a, $a=(a_1,b_1,c_1),$ $ b=(a_2,b_2,c_2)$, $c_1c_2\neq 0$. Then the condition 
$Q_{abx}=0$ implies $2c_1c_2r=- (a_2c_1+c_2a_1)p.$
If there exists $c=(a_3,b_3,c_3)\in C$ such that $c_3a_2\neq  a_3c_2$  then
 $0=Q_{acy}=-c_1p(c_3a_2-a_3c_2)/c_2$ gives $p=0, r=0$ which proves the lemma. 
If  $a_3c_2=a_2c_3 $  we can suppose that $b_3c_2\neq c_3b_2,\;$ and the condition 
$ 0=Q_{acx}=c_1p(c_2b_3-b_2c_3)/c_2$ gives  $r=p=0$.
\smallskip

In the case 2b we get  $a=(a_1,b_1,0), b=(a_2,b_2,0), a_1a_2\neq 0$, and hence 
$Q_{abx}=3a_1a_2\neq 0$.

\smallskip\smallskip
3). Reducible case, $Q= z(x^2+py^2+qz^2+rxz+syz)$. Subcases:

\smallskip \hskip .5 cm 3a). The curve $C$ is not in ${\S}^2\bigcap\{z=0\}.$
 
\smallskip
\hskip .5 cm 3b). The curve $C\subset {\S}^2\bigcap\{z=0\}.$

\smallskip Suppose 3a, $a=(a_1,b_1,c_1),$ $ b=(a_2,b_2,c_2)$, $c_1c_2\neq 0$. Then the condition 
$ Q_{aby}=0$ implies that $ c_1c_2s=- (b_2c_1+c_2b_1)p.  $
For any $c=(a_3,b_3,c_3)$ one gets $0= Q_{acy}={p(b_3c_2-c_3b_2)c_1}/{c_2}$ with 
 $b_3c_2\neq c_3b_2$ since $b_3c_2= c_3b_2$ gives $Q_{acx}= {(a_3c_2-c_3a_2)c_1}/{c_2}\neq 0$.
 Hence $s=p=0$.
                          
\smallskip Suppose 3b,  $a=(a_1,b_1,0),\; b=(a_2,b_2,0),\; c=(a_3,b_3,0), \; a_1a_2\neq 0,\;$\linebreak
$ b_1b_2\neq 0,\; a_2b_3\neq  a_3b_2.$
Then $0=Q_{abz}=pb_1b_2+a_1a_2,$ $p= -a_1a_2/(b_1b_2)$, $Q_{acz}=a_1(a_3b_2-b_3a_2)/b_2\neq 0,$
a contadiction and the proof is finished.

\smallskip

\section{Proof of the Theorem }

We begin with the following construction.
\smallskip
 
 Let $x\in \S^3$. Set
$$A_x = \{ (a,b)\in \S^3\times \S^3, a\neq b: u_{a,b} (x) = \sup_{y\in \S^3} u_{a,b}(y)\} ;$$
note that $A_x$ is a semi-analytic subset of $ \S^3\times \S^3$, and $(a,b) \in A_x $ implies 
$(b,a) \in A_x $. The semi-analycity of $A_x$ implies the sub-analycity of all the sets below in the proof.
In particular they verify  Whitney's stratification theorem [W] as was showed by Hironaka [H],
i.e. each such set $M$ is  stratified in a finite union of open $k$-dimensional smooth submanifolds, 
$k=0,1,...,m=\dim M$.

\smallskip
Let then  ${\frak C}^x$ for $x\in \R^4  \setminus \{0\}$ be the cubic form of the Taylor expansion of the 
function $u$ at the point $x$,  i.e., $D^3{\frak C}^x= D^3u(x) .$ Let us notice first that for any vector 
$e\in \R^4$ the function $u_e$ is a homogeneous order 1 and hence $x$ is a zero eigenvector of the
 quadratic form $({\frak C }^x_e)$.  We need  the following two simple properties of this form.

\medskip
{\bf Lemma 5.} {\it Let $(a,b) \in A_x $. Then $b$ is a zero eigenvector of the
quadratic form ${\frak C}^x_a$. }

\medskip{\em Proof.} From our assumptions it  follows that for any vector $e\in \R^4$ one has $u_{a,b,e}(x)=0$. Hence $({\frak C }^x_a)_{b,e}=0$. This implies that $b$ is a zero eigenvector of ${\frak C }^x_a$.

\medskip{\bf Lemma 6.} {\it Let   $a,x,b_1, b_2, b_3 \in \S^3$ with linearly independent $b_1, b_2, b_3$
 such that  $(a,b_1), (a,b_2),(a,b_3)\in A_x $. Then ${\frak C }^x_a=0$. }

\medskip{\em Proof.} By Lemma 5 the vectors $b_i$ are zero eigenvectors of the quadratic form 
${\frak C }^x_a$, i.e., it has the  zero eigenvalue with  multiplicity at least 3. Since ${\frak C }^x_a$
 should change the sign or be equal  zero the lemma follows.

\medskip Let now
 $$X:=\{x\in  \S^3 : \dim A_x \geq 3\}.$$
 Then $X\neq\emptyset$ since
$$\bigcup_{x\in\S^3}A_x= \S^3 \times \S^3,\;\dim( \S^3 \times \S^3) =6,$$ 
we denote $d\in [0,3]$ its dimension.
  
 \smallskip  Let $\Gamma =\cup_{x\in X }A_x$  then $
\dim(\S^3 \times \S^3\setminus \Gamma)\le 5, \dim(\Gamma) =6$.  

\smallskip
We have four possibilities for $d$, namely, $d=0,1,2$ or 3.

\bigskip

1. Let $d=0$. Then $\dim A_y=6$ for some $y\in X$, and 
  $$\dim((\S^3 \times \{e\} )\bigcap A_y)\ge 3$$  
for $e\in \S^3$.

\smallskip In this case one can find  linearly  independent vectors $e_1,...,e_4$, $e_i\neq y$, 
such that $(e,e_i)\in A_y$. Applying Lemma 1 to the function $u_e$ we get the proof.

\smallskip 2.  Let $d=1$. Then we can suppose without loss that $\dim A_y= 5$ for any  $y\in X$ and
$$\dim((\S^3 \times \{e\} )\bigcap A_y)\ge 2, \;\dim(( \{e\} \times \S^3)\bigcap A_y)\ge 2$$  
thus
$$E_1\times E_2\subset A_y$$
  $ E_1, E_2\subset \S^3, \dim(E_1)=\dim(E_2)=2.$    

\smallskip   Denote the set of all $y\in \S^3$ satisfying $E_1\times E_2\subset A_y$ by $Y$. 
Let $y\in Y$, $a\in E_1$. Then By Lemma 6 ${\frak C}^y_a=0$.  Since $E_1 $ is a  2-dimensional
 set the cubic form  ${\frak C}^y$ depends at most on one coordinate. Since its
derivative change sign it follows that  ${\frak C}^y=0$. Thus if $Y_1$ is a connected component
of $Y$ then $D^2u$ is constant on $Y_1$.   On the other hand since $Y$ is a real analytic set it contains only
finite number of connected components, $Y_1,...,Y_n$. At each $Y_i$ function $u$ has a fixed
Hessian. Therefore there is at least one $Y_j$ such that for $y\in Y_j$ the set $A_y$ is 6-dimensional
and one returns to the previous case.
 
\smallskip 3.  Let $d=2$. We  suppose without loss that $\dim A_y= 4$ for any  $y\in X.$ 
For a connected component $A,\; \dim A=4$ of $A_y$ let $d_1=d_1(A),\;d_2=d_2(A)$ 
be the dimensions of the  projections of $A$ to the first and the second factor  in the product 
$\S^3 \times \S^3$ respectively. By symmetry one can suppose $d_1\ge d_2\ge 1.$
 Since $d_1+d_2\ge\dim A= 4$ we have the following possibilities:

\smallskip
3a). $d_1=3,d_2=1;$

\smallskip
3b). $d_1=2,d_2=2;$

\smallskip
3c). $d_1=3,d_2=2;$

\smallskip
3d). $d_1=d_2=3.$

\smallskip Since in the cases 3a and 3b one has $d_1+d_2=\dim A$, the manifold $A$ itself is a product and we return
to the cases 1 and 2 respectively.

\smallskip
Suppose 3c or 3d and let $Z\subset \S^3$ be the image of the first projection of $A_x, \;\dim Z=3.$
Then for any $x\in Z$ there is a curve $\gamma_x \subset \S^3$ verifying the following condition:
$$\forall a\in \gamma_x, \; a\times D(a)\subset  A_x$$
for a 1- or 2-dimensional set $D(a)\subset \S^3 .$

Let $y\in Z$, and let $a,a' \in \gamma_y , a\neq a'$. Then By Lemma 5 ${\frak C}^y_a=0$,  
${\frak C}^y_{a'}=0$  and hence  ${\frak C}^y$ does not depend  on the coordinates parallel 
to $a$ and $a'$. Thus the cubic form  ${\frak C}^y$ depends at most on two coordinates.
Thus for any $e\in  \S^3$ the rank of the gradient map  $\nabla {\frak C}^y_e \rightarrow \R^4$
is at most 2 at the point $y\in Z$. Therefore  since $u_e$ is a homogeneous order one function the rank 
of the gradient map $\nabla_x u_e: \S^3 \rightarrow \R^4$ is at most 2 at any point $y\in Z$.
For   an affine hyperplane  $L\subset \R^4, 0\notin L$ let $Z'$ be the spherical projection of $Z$ on $L$, and let
 $s =u_e|_L$. Since $u_e$ is a homogeneous order one function the gradient map of $u_e(x)$ depends only 
on the spherical coordinate of $x$ it follows that $\det D^2s =0$ on $Z'$.
Since $s$ is a real analytic function and $Z'$ is a 3-dimensional we get $\det D^2s =0$ on the whole plane $L$ and thus  by Lemma 2 $u_e$ is linear.

 \smallskip 4.  Let $d=3$.   We  suppose without loss that $\dim A_y= 3$ for any  $y\in X.$ For a connected
component $A,\dim A=3$ of $A_y$ let $d_1\ge d_2$ be as before, $d_1+d_2\ge 3$. One has the following possibilities:

\smallskip
4a). $d_1=2,d_2=1;$

\smallskip
4b). $d_1=d_2=2;$

\smallskip
4c). $d_1=3,d_2=0;$

\smallskip
4d). $d_1=3,d_2=1;$ 

\smallskip
4e). $d_1=3, d_2=2;$

\smallskip
4f). $d_1=d_2=3.$

\smallskip
In the case 4a one has $A_x=E_1\times C_2, \dim E_1=2, \dim C_2=1$ and the proof above for
$A_x=E_1\times E_2,\; \dim E_1=\dim E_2=2$ remains valid.

\smallskip
In the case 4c  one has $A_x=\S^3\times \{a\}$ and we return to the case 1.

\smallskip
Suppose then 4d, 4e or 4f, let  $Z_x: =pr_1(A_x)\subset \S^3,\; \dim Z_x=3$
Then for any $x\in X$ one gets:
 $$\forall a\in Z_x, \; a\times h(a)\in  A_x,$$
where $h(a)\in \S^3$.

 \smallskip   Let $y\in X$ and let $L=y^{\perp}\subset \R^4$. Since $u$ is a homogeneous order 2
 function ${\frak C}^y$ depends only on the
coordinates of $L$. Thus there exists a 2-dimensional set $E\subset \S^2 \subset L$ such that
${\frak C}^y_e$ is degenerate for any $e\in E$ and hence for any $e\in \S^2$. Thus by Lemma 3
the cubic form ${\frak C}^y$ depends only on 2 variables and we  finish the proof  as for $d=2.$

Assume finally 4b, and let $y\in X$.

Then by Lemma 4 the cubic form ${\frak C}^y$  depends only on 2
coordinates, which we denote by $z_1, z_2$; let $l$ be the linear span of $z_1, z_2$. 
Thus  $l$ is a zero eigenspace of  ${\frak C}^y_e$ for any $e\in \S^3$ . By our assumption 
 one finds $(a,b)\in A_y, \; b  \notin l$. Therefore the multiplicity of the zero eigenvalue 
of ${\frak C}^y_a$ is at least 3. Again, since its derivatives change sign it follows that  
${\frak C}^y=0$ and one finishes the proof as before. 
 
  \medskip \medskip\centerline{REFERENCES}

 \medskip

\noindent [A] A.D. Alexandroff, {\it Sur les th\'eor\`emes
d'unicite pour les surfaces ferm\'ees}, Dokl. Acad. Nauk 22
(1939), 99--102.

\medskip
 \noindent [ASS] S. N. Armstrong, L. Silvestre, C. K. Smart, {\em Partial regularity of solutions of fully nonlinear uniformly elliptic equations}, arXiv:1103.3677.

\medskip
 \noindent [C] L. Caffarelli,  {\it Interior a priory estimates for solutions
 of fully nonlinear equations}, Ann. Math. 130 (1989), 189--213.

\medskip
 \noindent [CC] L. Caffarelli, X. Cabre, {\it Fully Nonlinear Elliptic
Equations}, Amer. Math. Soc., Providence, R.I., 1995.

\medskip
 \noindent [CIL]  M.G. Crandall, H. Ishii, P-L. Lions, {\it User's
guide to viscosity solutions of second order partial differential
equations,} Bull. Amer. Math. Soc. (N.S.), 27(1) (1992), 1--67. 

 \medskip
 \noindent [CL] S.-S. Chern, R. Lashof, {\it On the total curvature of immersed manifolds, }
 Amer. J. Math. 79 (1957), 306-318.

 \medskip
 \noindent [GT] D. Gilbarg, N. Trudinger, {\it Elliptic Partial
Differential Equations of Second Order, 2nd ed.}, Springer-Verlag,
Berlin-Heidelberg-New York-Tokyo, 1983.

\medskip
 \noindent [HNY] Q. Han, N. Nadirashvili, Y. Yuan, {\it Linearity of
homogeneous order-one solutions to elliptic equations in dimension
three,} Comm. Pure Appl. Math. 56 (2003), 425--432.

\medskip
 \noindent [H] H. Hironaka, {\it Normal cones in analytic Whitney stratification, } Inst. Hautes
 \'Etudes Sci. Publ. Math. 36 (1969), 127-139.

\medskip
 \noindent [LO] H.B. Lawson, Jr., R. Osserman, {\it Non-existence, non-uniqueness and irregularity of solutions
to the minimal surface system.}  Acta Math. 139 (1977), 1--17.

 \medskip
 \noindent [M] J. S. Milne,  {\it Elliptic curves,} BookSurge Publishing, 2006.

\medskip
\noindent [NV1] N. Nadirashvili, S. Vl\u adu\c t, {\it
Nonclassical solutions of fully nonlinear elliptic equations,}
Geom. Func. An. 17 (2007), 1283--1296.

\medskip
\noindent [NV2] N. Nadirashvili, S. Vl\u adu\c t, {\it Singular Viscosity Solutions to Fully Nonlinear
Elliptic Equations},  J. Math. Pures Appl., 89 (2008), 107-113.

\medskip
\noindent [NV3] N. Nadirashvili, S. Vl\u adu\c t, {\it Octonions and Singular Solutions of Hessian  Elliptic
Equations}, Geom. Func. An. 21 (2011), 483-498.

\medskip
\noindent [NV4] N. Nadirashvili, S. Vl\u adu\c t, {\it Singular solutions of Hessian fully nonlinear elliptic
equations}, Adv. Math., 228 (2011), 1718-1741.

\medskip
\noindent [NTV] N. Nadirashvili, V. Tkachev, S. Vl\u adu\c t, {\it Non-classical Solution to Hessian   
Equation   from Cartan Isoparametric Cubic}, arXiv: 1111.0329

\medskip
\noindent [S] S. Sternberg, {\it Lectures on Differential Geometry, } Prentice-Hall, 1964.

\medskip
\noindent [T] N. Trudinger, {\it H\"older gradient estimates for
fully nonlinear elliptic equations,} Proc. Roy. Soc. Edinburgh
Sect. A 108 (1988), 57--65.

   \medskip
 \noindent [W] H. Whitney, {\it  Local properties of analytic varieties, } Differential and Combinatorial Topology, Prrinceton Univ. Press (1965), 205-244.

\end{document}